\documentclass[12pt]{article}
\usepackage{amsfonts,amsmath,amssymb,amscd,amsthm}
 \def\R{{\mathbb R}}
 \def\Q{{\mathbb Q}}

\long\def\comment#1\endcomment{}

 \def\R{{\mathbb R}}
 \def\Q{{\mathbb Q}}

\long\def\comment#1\endcomment{}

\begin{document}
\centerline{\uppercase{\bf Solvability of cubic and quartic equations
using one radical}}
\footnote {This paper is prepared under the supervision of Arkadiy Skopenkov and is submitted to the Moscow Mathematical Conference for High-School Students. Readers are invited to send their remarks and reports on this paper to mmks@mccme.ru. }

\bigskip

\centerline{\bf Danil Akhtyamov, Ilya Bogdanov}

\bigskip

Theorem 1.1 from [CK] is  similar to Theorem 2 from this note. And Theorem 1.2 from [CK] is similar to Theorem 3 from this note. But Theorems 2 and 3 do not follow from Theorems, which proved in [CK]. And Theorems of Chu and Kang do not follow from Theorems 2 and 3.

\comment
The note we present shorter proofs of the `only if' part of  Theorems 1 and 2 below essentially proved in [K, CK], cf [A].
\footnote{Clearly, Theorem 1 follows from [K].
The first condition of Theorem 2 follows from [CK] because for the equation $T^4-uT^2-vT-w=0$
the {\it cubic resolution} $4(u+2\alpha)(w+\alpha^2)-v^2$ is obtained from $P(X):=64X^3-32uX^2+(4u^2+16w)X-v^2$
by a linear change of variable $X:=\frac{2\alpha+u}{4}$:
$$P(\frac{2\alpha+u}{4})=
((2\alpha+u)^2-2u(2\alpha+u)+(u^2+4w))(2\alpha+u)-v^2=
4(w+\alpha^2)(u+2\alpha)-v^2.$$
}

\endcomment

\section{Formulations}

We say that a polynomial is {\it $1$-solvable} if it has a root $x\in\Q[r]$, where $r\in\R$ and $r^n\in\Q$ for some positive integer $n$. A `polynomial with rational coefficients' is referred to merely as {\em polynomial}.

\smallskip

{\bf Theorem 1. }{\it Assume that a polynomial of degree~$n$ is irreducible over~$\Q$ and 1-solvable. Then this polynomial has a root of the form $A(r)$, where $A$ is a polynomial, and a number $r\in\R$ satisfies $r^n\in\Q$.}

\smallskip
{\bf Theorem 2.}
{\it For a cubic polynomial $x^3+px+q$, the following conditions are equivalent:

(1-solvability) it is 1-solvable;

($a+br+cr^2$) this equation has a root of the form $a+br+cr^2$, where $r\in\R$ and $a,b,c,r^3\in\Q$;

($\sqrt{D_{pq}}\in\Q$) either it has a rational root, or $D_{pq} \ge 0$ and $\sqrt{D_{pq}}\in\Q$, where $D_{pq}=(\frac{p}{3})^3+(\frac{q}{2})^2 $}

{\bf Theorem 3.}
{\it For an irreducible  polynomial $x^4+px^2+qx+s$ , the following conditions are equivalent:

(1-solvability) it is 1-solvable;

($a+br+cr^2+dr^3$) this equation has a root of the form $a+br+cr^2+dr^3$, where $r\in\R$ and $a,b,c,d,r^4\in\Q$;

($\sqrt{\Gamma}\in\Q$) there exists $\alpha\in\Q$ such that $2\alpha>p$ and $q^2-4(p-2\alpha)(s-\alpha^2)=0$, and the number $\Gamma=16(\alpha^2-s)^2-(\alpha^2-s)(2\alpha+p)^2$ is a square of a rational number.}

\smallskip

\section{Proof of Theorem 1.}

\smallskip

Proof of Theorem 1 uses the following statement.

\smallskip

{\bf Statement 1. }Assume that $r\in\R\setminus\Q$ and $r^n\in\Q$ for some integer $n>1$. Take any $\beta\in\Q[r]$. Then there exists a positive integer~$k$ such that $\beta\in\Q[r^k]$ and $r^k\in\Q[\beta]$ (in other words, $\Q[r^k]=\Q[\beta]$).

\smallskip

{\it Proof. }The number~$r$ is a root of some nonzero polynomial with coefficients in~$\Q[\beta]$ (e.g., $x^n-r^n$. Choose such polynomial~$f(x)$ of the least possible degree~$k$.

Consider the g.c.d. of~$x^n-r^n$ and~$f$; it also has $r$ as a root, its coefficients lie in~$\Q[\beta]$, and its degree does not exceed~$k$; this means that this g.c.d. is~$f$ itself. So all the complex roots of~$f$ have the form~$r\varepsilon_n^i$. Then, by the Vieta Theorem, the absolute value of the constant term of~$f$ equals~$r^k$ for some $k\leq n-1$. Since this constant term is real, we obtain that $r^k\in\Q[\beta]$.

Now it remains to prove that $\beta\in\Q[r^k]$. Since $\beta\in\Q[r]$, we have
$$
\beta=b_0(r^k)+rb_1(r^k)+\ldots+r^{k-1}b_{k-1}(r^k)
$$
for some polynomials~$b_0,\ldots,b_{k-1}\in\Q[x]$. If not all polynomials~$b_1,\ldots,b_{k-1}$ are zeroes, then $r$ is a root of a nonzero polynomial $$(b_0(r^k)-\beta)+xb_1(r^k)+\ldots+x^{k-1}b_{k-1}(r^k)$$
whose degree is~$k$, and whose coefficients lie in~$\Q[\beta]$ (since $\beta,r^k\in\Q[\beta]$). This contradicts the choice of~$f(x)$. Thus we arrive at $b_1=\dots=b_{k-1}=0$, whence $\beta=b_0(r^k)\in\Q[r^k]$.

\smallskip

{\it Proof of Theorem 1. }Let the given polynomial has a root $y_0\in\Q[R]$ for some $R\in\R$ satisfying $R^D\in\Q$, where $D$ is positive integer number. By~Statement 1, we have $\Q[y_0]=\Q[R^k]$ for some positive integer $k$.

Denote $r=R^k$. Since $R^D\in\Q$, one may choose the minimal positive integer~$d$ such that $r^d\in\Q$; then $r,r^2,\dots,r^{d-1}\notin\Q$. Therefore, the polynomial $x^d-r^d$ is irreducible over~$\Q$ (since the  constant term of any its nontrivial unitary factor has an irrational absolute value~$r^t$, $0<t<d$).

Finally, the equality~$\Q[y_0]=\Q[r]$, combined with the dimension argument 
yield that any two irreducible (over~$\Q$) polynomials, one with root~$y_0$ and the other with root~$r$, have equal degrees. This shows that $n=d$, as required. QED.
\smallskip

\section{Proof of Theorem 2.}

\smallskip

\underline{\it (1-solvability) $\Rightarrow (a+br+cr^2)$.} If the given polynomial is reducible, then it has a rational root which has the required form. Otherwise the result follows directly from Theorem 1.

\smallskip

\underline{\it $(\sqrt{D_{pq}}\in\Q)\Rightarrow$ (1-solvability).} Set $r=\root3\of{-\frac{q}{2}+\sqrt{D_{pq}}}$. By Cardano's formula, the unique real root of the equation $x^3+px+q=0$ equals
$$
  r-\frac{p}{3r}=r-\frac{p}{3r^3}\cdot r^2=r-\frac{p}{3\bigl(-\frac q2+\sqrt{D_{pq}}\bigr)}\cdot r^2.
$$

\smallskip

Proof of the '{\it $(a+br+cr^2)\Rightarrow(\sqrt{D_{pq}}\in\Q)$}' part of Theorem 2 uses the following Lemmas and Statements:

{\bf Conjugation Lemma}.
{\it Let $a,b,c,r^3\in\Q$ , $r\not\in\Q$  and number $x_0:=a+br+cr^2$ is a root of polynomial $F$ with rational coefficients.
Then numbers $$ x_1=a+br\epsilon+cr^2\epsilon^2, \quad x_2=a+br\epsilon^2+cr^2\epsilon.$$
are also roots of $F$.}

\smallskip

\underline{\it $(a+br+cr^2)\Rightarrow(\sqrt{D_{pq}}\in\Q)$.} Suppose the contrary.
If $r\in\Q$ or $b=c=0$, then the equation has a rational root. In the remaining case, denote $\varepsilon=\frac{i\sqrt{3}-1}{2}$. Each of the numbers $x_1$, $x_2$, and~$x_3$ defined in Conjugation Lemma is a root of our equation. By $D_{pq} \ne 0$, these three roots are distinct. Therefore, $x_1$, $x_2$, and $x_3$ are all the roots of our equation. Now,  we have
\begin{multline*}
   -108D_{pq}=(x_2-x_3)^2(x_1-x_3)^2(x_1-x_2)^2\\
   =\bigl(br(\varepsilon-\varepsilon^2)+cr^2(\varepsilon^2-\varepsilon)\bigr)^2
     \bigl(br(1-\varepsilon)+cr^2(1-\varepsilon^2)\bigr)^2
     \bigl(br(1-\varepsilon^2)+cr^2(1-\varepsilon)\bigr)^2\\
   =\varepsilon^2(1-\varepsilon)^6(br-cr^2)^2
     \bigl(br+cr^2(1+\varepsilon)\bigr)^2
     \bigl(br(1+\varepsilon)+cr^2\bigr)^2.
\end{multline*}
Since $(1+\varepsilon)(1+\varepsilon^2)=(-\varepsilon)(-\varepsilon^2)=1$
and $(\varepsilon-1)^3=3\varepsilon-3\varepsilon^2=3\sqrt{3}\,i$, we obtain
\begin{multline*}
   -108D_{pq}
   =-27\varepsilon^2(1+\varepsilon)^2(br-cr^2)^2
     \bigl(br(1+\varepsilon^2)+cr^2\bigr)^2
     \bigl(br(1+\varepsilon)+cr^2\bigr)^2\\
   =-27(\varepsilon+\varepsilon^2)^2(br-cr^2)^2
     \bigl(b^2r^2+br\cdot cr^2+c^2r^4\bigr)^2
   =-27\bigl((br)^3-(cr^2)^3\bigr)^2.
\end{multline*}
This yields the required result. QED.

\section{Proof of Theorem 3.}

\smallskip

\underline{\it (1-solvability) $\Rightarrow$ ($a+br+cr^2+dr^3$) .} It follows from Theorem 1.

\smallskip

Proof of the'{\it ($\sqrt{\Gamma}\in\Q$) $\Rightarrow$ ($a+br+cr^2+dr^3$)'} part of
Theorem 3 uses the following well-known Statement:

\smallskip

{\bf Statement(Ferrari's Formula). }{\it Assume that $F(x)=x^4+px^2+qx+s$ has rational coefficients, $\alpha\in\R$, $q^2-(p-2\alpha)(s-\alpha^2)=0 $, $2\alpha-p>0$ and $$-2\alpha-p+4\sqrt{\alpha^2-s} >0$$.
Let $$x_\pm:=\pm\sqrt{2\alpha-p}+\sqrt{-2\alpha-p+4\sqrt{\alpha^2-s}}$$. Then either $F(x_{+})=0$ or $F(x_{-})=0$.}

\smallskip

{\it Proof.}   We have $q=\pm\sqrt{(2\alpha-p)(\alpha^2-s)}$. Then

$$x^4+px^2+qx+s=(x^2+\alpha)^2-(x\sqrt{2\alpha-p}\pm\sqrt{\alpha^2-s})^2=$$
$$=(x^2+\alpha+x\sqrt{2\alpha-p}\pm\sqrt{\alpha^2-s})(x^2+\alpha-x\sqrt{2\alpha-p}\mp\sqrt{\alpha^2-s}). \quad\text{QED.}$$

\smallskip

\underline{\it ($\sqrt{\Gamma}\in\Q$) $\Rightarrow$ (1-solvability).} By $\alpha^2-s>0$ and $\Gamma \ge 0$ we have $16(\alpha^2-s)-(2\alpha+p)^2 \ge 0$.
 Let $r$ be real number such that
$$2r^2=4\sqrt{\alpha^2-s}+\sqrt{16(\alpha^2-s)-(2\alpha+p)^2}.$$
 Then $r^4\in\Q$, because $$4r^4=16(\alpha^2-s)+8\sqrt{16(\alpha^2-s)^2-(\alpha^2-s)(2\alpha+p)^2}+$$
$$+16(\alpha^2-s)-(2\alpha+p)^2=32(\alpha^2-s)+8\sqrt{\Gamma}-(2\alpha+p)^2\in\Q.$$
Let $F(t):=t^2-4t\sqrt{\alpha^2-s}+\frac{(2\alpha+p)^2}{4}$. Then $F(r^2)=0$, i.e. $r^4-4r^2\sqrt{\alpha^2-s}+\frac{(2\alpha+p)^2}{4}=0$. Then by $r^4\in\Q$ we have $r^2\sqrt{\alpha^2-s}\in\Q$.
Then $4\sqrt{\alpha^2-s}=r^2+\frac{(2\alpha+p)^2}{4r^2}$.
Then $$\sqrt{-2\alpha-p+4\sqrt{\alpha^2-s}}=\sqrt{-2\alpha-p+r^2+\frac{(2\alpha+p)^2}{4r^2}}=$$.

$$=\sqrt{(r+\frac{-2\alpha-p}{2r})^2}=|r+\frac{-2\alpha-p}{2r}|=|r+\frac{r^3(-2\alpha-p)}{2r^4}|.$$
We have $\pm\sqrt{(2\alpha-p)(\alpha^2-s)}=q\in\Q$ and  $r^2\sqrt{\alpha^2-s}\in\Q$. Then $r^2\sqrt{2\alpha-p}\in\Q$. Then $x_\pm$  are numbers of the form $br+cr^2+dr^3$, where $b,c,d,r^4\in\Q, r\in\R$. QED.

\smallskip

Proof of the `{\it ($a+br+cr^2+dr^3)$ $\Rightarrow$ ($\sqrt{\Gamma}\in\Q$)}' part of Theorem 3 uses the following well-known Lemmas and Statements, some of which we prove for completeness.

\smallskip

\comment
{\bf Linear Independence Lemma.}
{\it If $r\in\R$, $r^2\not\in\Q$, $r^4, w, u, v, t\in\Q$  and $wr^3+ur^2+vr+t=0$, then $w=u=v=t=0$.}
\endcomment
\comment

\smallskip
{\it Proof}. Suppose the contrary.
 Let divide $x^4-r^4$ to $wx^3+ux^2+vx+s$:
$x^4-r^4=(wx^3+ux^2+vx+s)A(x)+B(x)$.
Then $B(r)=0$.
There are 2 cases:

1)$B(x)=0$,

2)$B(x) \ne 0$

Consider the first case. Then $x^4-r^4$ reducible over $Q$. A contradiction.

Consider the second case. Let divide $x^4-r^4$ to $B(x)$: $x^4-r^4=B(x)C(x)+D(x)$.
Then $D(x)=0$, because $D(r)=0$ and $D(x)$ is the first degree polynomial or a constant. Then $x^4-r^4$ reducible over $Q$. A contradiction. QED.

\endcomment

\smallskip
{\bf Conjugation Lemma}.
{\it Let $a,b,c,d,r^4\in\Q$ , $r^2\not\in\Q$  and number $x_0:=a+br+cr^2+dr^3$ is a root of polynomial $F$ with rational coefficients.
Then numbers
$$x_1:=a-br+cr^2-dr^3,\quad x_2:=a+bri-cr^2-dr^3i,\quad x_3:=a-bri-cr^2+dr^3i$$
are also roots of $F$.}

\smallskip
{\it Proof}. Let $F(a+bx+cx^2+dx^3)=B(x)(x^4-r^4)+wx^3+ux^2+vx+t$ be a division with a remainder.
Then $wr^3+ur^2+vr+t=0$.
By the Linear Independence Lemma $w=u=v=t=0$.
Then $F(a+bx+cx^2+dx^3)=B(x)(x^4-r^4)$.
Then $F(a-br+cr^2-dr^3)=F(a+bri-cr^2-dr^3i)=F(a-bri-cr^2+dr^3i)=0$. QED.

\smallskip
{\bf Statement 2}.{\it Let $r\in\R$, $b,c,d,r^4 \in\Q$, $r^2\not\in\Q$,  $b^2+d^2\ne 0$. Then

(a) $br-dr^3 \ne 0$ and $br+dr^3 \ne 0.$

(b) numbers
$x_0,x_2,x_1,x_3$ from the Conjugation Lemma.
are pairwise different.}

\smallskip
{\it Proof of (a)}. Assume the contrary.
If $d=0$, then $b=0$, which is a contradiction.
If $d\ne0$, then either $r^2=b/d\in\Q$ or $r^2=-b/d\in\Q$, which is a contradiction. QED.

\smallskip
{\it Proof of (b)}. Suppose the contrary.
Then either $x_0=x_2$ or $x_1=x_3$.
Then either $br-dr^3=0$ or $br+dr^3=0$, which contradicts to (a). QED.

\smallskip
{\bf Statement 3.}{ \it Let $x_0, x_1, x_2, x_3$ be the roots of polynomial $x^4+px^2+qx+s$ (not necessarily coinciding with the numbers from the Conjugation Lemma). Let $\alpha=\frac{x_0x_2+x_1x_3}{2}$. Then:

(a) $4(\alpha^2-s)=(x_0x_2-x_1x_3)^2$.

(b) $p-2\alpha=(x_0+x_2)(x_1+x_3).$

(c) $q^2-4(p-2\alpha)(s-\alpha^2)=0$ .

(d) $\Gamma=-(x_0x_2-x_1x_3)^2(x_1-x_3)^2(x_0-x_2)^2.$
 }

\smallskip

{\it Proof of (a)}.  By Vieta Theorem: $$4(\alpha^2-s)=(x_0x_2+x_1x_3)^2-4x_1x_2x_3x_4=(x_0x_2-x_1x_3)^2.\quad\text{QED.}$$

{\it Proof of (b)}.  By Vieta Theorem: $$p-2\alpha=x_0x_1+x_0x_3+x_1x_2+x_2x_3=(x_0+x_2)(x_1+x_3).\quad\text{QED.}$$

{\it Proof of (c)}. Since $x_0+x_1+x_2+x_3=0$, by 1(a) and 1(b) we have: $$q^2-4(p-2\alpha)(s-\alpha^2)=(x_0x_1x_2+x_0x_1x_3+x_0x_2x_3+x_1x_2x_3)^2+4(x_0+x_2)(x_1+x_3)\frac{(x_0x_2-x_1x_3)^2}{4}=$$
$$=(x_0x_2(x_1+x_3)+x_1x_3(x_0+x_2))^2+(x_0+x_2)(x_1+x_3)(x_0x_2-x_1x_3)^2=$$
$$=(x_0x_2-x_1x_3)^2(x_0+x_2)^2-(x_0+x_2)^2(x_0x_2-x_1x_3)^2=0.\quad\text{QED.}$$

{\it Proof of (d)}. Since $x_0+x_1+x_2+x_3=0$, by Statement 3(b) and 3(c): $$4\Gamma=4(x_0x_2-x_1x_3)^4-(x_0x_2-x_1x_3)^2((x_0+x_2)(x_1+x_3)+2(x_0x_2+x_1x_3))=$$
$$=(x_0x_2-x_1x_3)^2(2x_0x_2-2x_1x_3-2x_0x_2-2x_1x_3+(x_1+x_3)^2)\times$$
$$\times(2x_0x_2-2x_1x_3+2x_0x_2+2x_1x_3-(x_0+x_2)^2)=$$
$$=(x_0x_2-x_1x_3)^2((x_1+x_3)^2-4x_1x_3)(4x_0x_2-(x_0+x_2)^2)=$$
$$=-(x_0x_2-x_1x_3)^2(x_1-x_3)^2(x_0-x_2)^2$$

{\it \underline{\it ($a+br+cr^2+dr^3)$ $\Rightarrow$ ($\sqrt{\Gamma}\in\Q$)}}. By the Conjugation Lemma the numbers $x_0,x_1,x_2,x_3$ from the Lemma
are roots of $F$.
Then by the Statement 2(b) and Vieta Theorem $-4a=0$.
Let $\alpha:=\frac{x_0x_2+x_1x_3}{2}$.
Then polynomial $q^2-4(p-2x)(s-x^2)$ has a rational root $x=\alpha$, because
$$2\alpha=x_0x_2+x_1x_3=(cr^2)^2-(br+dr^3)^2+(cr^2)^2-(dr^3i-bri)^2=$$
$$=2c^2r^4-r^2(b+dr^2)^2+r^2(b-dr^2)^2=2c^2r^4-4r^2bdr^2\in\Q.$$

By 3(b) $p-2\alpha=(x_0+x_2)(x_1+x_3)=-(x_0+x_2)^2<0.$

By 3(d): $$4\Gamma=-(x_0x_2-x_1x_3)^2(x_1-x_3)^2(x_0-x_2)^2=$$
$$-(c^2r^4-(br+dr^3)^2-c^2r^4-(br-dr^3)^2)^2(2br+2dr^3)^2(2bri-2dr^3i^2)^2=$$
$$=16r^4(-2b^2r^2-2d^2r^6)^2(b+dr^2)^2(b-dr^2)^2=$$
$$=64r^8(b^2+d^2r^4)^2(b^2-d^2r^4)^2=64r^8(b^4-d^4r^8)^2.\quad\text{QED.}$$

\bigskip

\newpage

\centerline{\bf References.}

\bigskip

[ABGSSZ] Akhtyamov, D., Bogdanov, I., Glebov, A., Skopenov, A., Strelstova, E., Zykin, A. (2015)
Solving equations using one radical. http://www.turgor.ru/lktg/2015/4/index.htm

[CK] Chu, H., Kang, M. (2002). Quartic fields and radical extensions. J. Symbolic Computation, 34, 83-89.

[K] Kang, M. (2000). Cubic fields and radical extensions. Am. Math. Monthly, 107, 254-256.

[S] Skopenkov, A. (2014). Some more proofs from the Book: solvability and insolvability of equations in radicals. http://arxiv.org/abs/0804.4357v6.

\end{document}